\documentclass[11pt]{amsart}
\usepackage{amssymb,amsmath,amsthm}
\usepackage[mathscr]{euscript}
\usepackage{enumerate, verbatim, url}
\newtheorem{theorem}{Theorem}

\newtheorem{conjecture}{\bf Conjecture}
\newtheorem{proposition}[theorem]{Proposition}

\theoremstyle{remark}

\newtheorem*{remark}{Remark}

\numberwithin{theorem}{section} \numberwithin{equation}{section}

\newcommand{\Fqt}{\mathbb{F}_q(t)}

\newcommand{\calO}{\mathcal{O}}

\newcommand{\mfa}{\mathfrak{a}}
\newcommand{\mfb}{\mathfrak{b}}
\newcommand{\vol}{\textnormal{vol}}
\newcommand{\Vol}{\textnormal{Vol}}

\newcommand{\R}{\mathbb{R}}
\newcommand{\C}{\mathbb{C}}
\newcommand{\F}{\mathbb{F}}
\newcommand{\Cl}{{\text {\rm Cl}}}
\newcommand{\Tr}{{\text {\rm Tr}}}

\newcommand{\Q}{\mathbb{Q}}

\newcommand{\Z}{\mathbb{Z}}
\newcommand{\N}{\mathbb{N}}
\newcommand{\PP}{\mathbb{P}}
\newcommand{\SL}{{\text {\rm SL}}}
\newcommand{\GL}{{\text {\rm GL}}}
\newcommand{\textmod}{{\text {\rm mod}}}

\newcommand{\calP}{\mathcal{P}}

\newcommand{\Disc}{{\text {\rm Disc}}}
\newcommand{\Stab}{{\text {\rm Stab}}}

\newcommand{\calR}{\mathcal{R}}

\newcommand{\sump}{\sideset{}{'}\sum}

\newcommand{\irr}{\textnormal{irr}}
\newcommand{\calF}{\mathcal{F}}
\newcommand{\Virr}{V^{\pm}_{\irr}}
\newcommand{\Res}{\textnormal{Res}}

\newcommand{\Aut}{\textnormal{Aut}}
\def\HH{\mathbb{H}}

\begin{document}
\title[Four Perspectives on Secondary Terms  in the Davenport-Heilbronn Theorems]
{Four Perspectives on Secondary Terms in the Davenport-Heilbronn Theorems}

\author{Frank Thorne}
\address{Department of Mathematics, University of South Carolina,
1523 Greene Street, Columbia, SC 29208}
\email{thorne@math.sc.edu}

\begin{abstract}
This paper is an expanded version of the author's lecture at the Integers Conference 2011. We discuss
the secondary terms in the Davenport-Heilbronn theorems on cubic fields
and $3$-torsion in class groups of quadratic fields. Such secondary terms had been conjectured by
Datskovsky-Wright and Roberts, and proofs of these or closely related secondary terms were
obtained independently by Bhargava, Shankar, and Tsimerman \cite{BST}, Hough \cite{H}, Zhao \cite{Z},
and Taniguchi and the author \cite{TT}.

In this paper we discuss the history of the problem and highlight the diverse methods used in \cite{BST, H, Z, TT} to address it.
\end{abstract}

\maketitle
\section{Introduction}
This paper concerns the following two theorems.

\begin{theorem}\label{thm_rc}
Let $N_3^{\pm}(X)$ count the number of cubic fields $K$ with $0 < \pm \Disc(K) < X$. We have
\begin{equation}\label{conj_rc}
N_3^{\pm} (X) = C^{\pm} \frac{1}{12 \zeta(3)} X + K^{\pm} 
\frac{4 \zeta(1/3)}{5 \Gamma(2/3)^3 \zeta(5/3)} X^{5/6} + o(X^{5/6}),
\end{equation}
where $C^- = 3$, $C^+ = 1$, $K^- = \sqrt{3}$, and
$K^+ = 1$. 
\end{theorem}

\begin{theorem}\label{thm_rc_torsion}
For any quadratic field
with discriminant $D$, let $\Cl_3(D)$ denote the 3-torsion subgroup of the ideal class group
$\Cl(\Q(\sqrt{D}))$,
we have
\begin{equation}\label{eqn_torsion1}
\sum_{0 < \pm D < X} \# \Cl_3(D) = \frac{3 + C^{\pm}}{\pi^2} X + K^{\pm}
\frac{8 \zeta(1/3)}{5 \Gamma(2/3)^3} \prod_p \Bigg(1 - \frac{p^{1/3} + 1}{p (p + 1)} \Bigg)
 X^{5/6} + o(X^{5/6}),
\end{equation}
and
where the sum ranges over fundamental discriminants $D$, the product is over all primes, and the constants are as before.
\end{theorem}

The main terms are due to Davenport and Heilbronn \cite{DH}, and the secondary term
in Theorem \ref{conj_rc} was conjectured by Roberts \cite{R}
and implicitly by Datskovsky and Wright \cite{DW3}. 

The secondary terms were expected to be difficult to prove. However, progress has recently made in
four independent works, using a variety of methods. Both results above have been proved
by Bhargava, Shankar, and Tsimerman \cite{BST}, using the geometry of numbers, 
and by Taniguchi and the present author \cite{TT}, using Shintani zeta functions. Hough \cite{H} has proved
a variation of Theorem \ref{thm_rc_torsion} by studying the distribution of Heegner points in the upper
complex plane, and Zhao \cite{Z} has obtained a variation of Theorem \ref{thm_rc} for function fields,
using algebraic geometry.\footnote{A proof of Theorem
\ref{thm_rc_torsion} does not appear in \cite{BST}, but the authors have shown me a proof. The error terms
in \cite{H} and \cite{Z} are larger than the secondary terms at present, but both of these approaches naturally 
{\itshape explain} the secondary term, and both authors are currently working on refining their methods.
}

The author's lecture at the 2011 Integers Conference explained (briefly) how
each of these four approaches all shed light on these secondary terms, and this paper is an expanded version
of our lecture. We begin with some background on counting fields and on counting cubic fields in particular. The
reader might skip to Section \ref{sec_four_approaches} for our discussion of the secondary
terms.

\section{Counting Fields in General}
How does one count number fields? This question has been addressed in many excellent
expository accounts, such as that by Cohen, Diaz y Diaz, and Olivier \cite{CDO}, or Chapter 6
of Bhargava's ICM proceedings article \cite{B_icm}. Therefore, our overview will be brief. To a large extent
we concentrate on the ``trivial'' problem of counting quadratic fields, as some observations on this problem
anticipate the methods developed to count cubic fields with good error terms.

The two most important theorems in the subject were proved by
Minkowski and Hermite:
\begin{theorem}[Hermite]
There are only finitely many number fields of bounded discriminant.
\end{theorem}
\begin{theorem}[Minkowski]
The discriminant of a number field $K$ of degree $n$ satisfies
\begin{equation}
|\Disc(K)| \geq \bigg( \frac{n^n}{n!} \bigg)^2 \bigg( \frac{\pi}{4} \bigg)^n.
\end{equation}
\end{theorem}
If we write $n = r + 2s$, where $r$ is the number of real embeddings of $K$, and $s$ is the
number of pairs of complex embeddings, we may replace 
$(\pi/4)^n$ with $(\pi/4)^{2s}$. Moreover, by 
Stirling's formula $n! \sim \sqrt{2 \pi n} \big( \frac{n}{e} \big)^n$ and
therefore
\begin{equation}
|\Disc(K)| \geq \Big( e^2 \cdot \frac{\pi}{4} + o(1) \Big)^n = \Big(5.803\dots + o(1)\Big)^n.
\end{equation}
The constant $\frac{e^2 \pi}{4}$ can be improved; see Odlyzko \cite{O}.
Conversely, this bound is sharp apart from the constant; Golod and Shafarevich \cite{GS} proved
the existence of an {\itshape infinite class field tower} of fields $K$ 
for which $|\Disc(K)|^{1/n}$ is constant.\footnote{Martinet \cite{M} gave the
example of $F = \Q(\zeta_{11} + \zeta_{11}^{-1}, \sqrt{-46})$, which he proved
to have an infinite 2-class field tower, for which $|\Disc(K)|^{1/n} = 92.2\cdots$ for
each field $K$. We refer to Lemmermeyer \cite{L} for a discussion of related issues, along with
a very thorough bibliography.}

To prove these theorems\footnote{See \cite{N}, Ch. III.2 for complete proofs.}, 
use the $r$ embeddings $K \hookrightarrow \R$ and the $s$ pairs of embeddings
$K \hookrightarrow \C$ to embed the ring of integers of $K$ as a lattice in $n$-dimensional
space. As proved by Minkowski, any convex body centered around the origin, whose volume is greater than
$2^n |\Disc(K)|$, must contain a lattice
point other than zero. However, this point must have norm at least 1, and geometric considerations allow one to conclude
Minkowski's lower bound.

Hermite's theorem can be proved in a similar way. We again use Minkowski's convex body theorem to find an element
$\alpha \in \calO_K$, the sizes of whose complex embeddings are all small, and which is guaranteed to generate $K$ over $\mathbb{Q}$.
The minimal polynomial of $\alpha$ is equal to $\prod_{\sigma} (X - \sigma(\alpha))$, where $\sigma$ ranges over all real
and complex embeddings of $K$, and the bounds obtained from Minkowski's theorem yield bounds on the coefficients
of the minimal polynomial of $\alpha$. For a fixed discriminant bound, there are only finitely many such polynomials, and thus
only finitely many possibilities for $K$.

These classical theorems suggest that it is natural to count fields by degree. In what follows, let $N_n(X)$
be the number of fields $K$ of degree $n$ with $|\Disc(K)| < X$. (Results are also known when the sign of the
discriminant, or more generally the number of real embeddings, is specified.)

We have the following results (excluding $n = 3$):
\\
\\
{\itshape n = 1.} There is only $\mathbb{Q}$.
\\
\\
{\itshape n = 2.} We have the equality of Dirichlet series
\begin{equation}\label{eqn_n2}
\sum_{[K : \Q] = 2} |\Disc(K)|^{-s} = -1 + 
 \Big( 1 - 2^{-s} + 2 \cdot 4^{-s}\Big) \frac{\zeta(s)}{\zeta(2s)},
\end{equation}
and it follows that 
\begin{equation}\label{eqn_n2_2}
N_2(X) = \frac{6}{\pi^2} X + O(\sqrt{X}).
\end{equation}

The typical proof of
\eqref{eqn_n2_2} estimates, for each squarefree odd $q$, the number of integers $n$
divisible by $q^2$ satisfying the 2-adic conditions implied by \eqref{eqn_n2_2}, and then uses inclusion-exclusion
to sieve for squarefreeness. The details are commonly left as an exercise for beginning analytic number theory
students. However, it was not until the recent work of Belabas, Bhargava, and Pomerance \cite{BBP} that
the analogous method was developed to yield power-saving error terms for $N_3(X)$.

How does one estimate the number of integers $n$ with $0 < n < X$ and $q^2 | n$? (For simplicity
of explanation, we ignore the 2-adic conditions.) One way is to interpret it as a lattice-point counting
problem in a one-dimensional vector space, and estimate the number of such $n$ by $\frac{1}{q^2}$
times the length of the interval $(0, X)$. This is of course trivial, but its generalization to the cubic case,
where one counts $\GL_2(\Z)$-orbits on a four-dimensional lattice, satisfying local
conditions $(\textmod \ q^2)$ and a global condition $0 < \pm \Disc(f) < X$, is not. (The higher
dimensional analogues are still less so!)

These integers may also be counted using analytic number theory. By Perron's formula, we have
\begin{equation}
\sum_{\substack{0 < n < X \\ q^2 | n}} 1 = \int_{2 - i \infty}^{2 + i \infty} q^{-2s} \zeta(s) X^s \frac{ds}{s},
\end{equation}
and these sums may be estimated by shifting the contour and using the functional equation of the zeta
function. This functional equation relies on Poisson summation, and so again relies on the question's
interpretation as a lattice point counting problem.

Although even the most diehard fan of contour integration would likely prefer the trivial proof,
this analytic method also generalizes usefully to the cubic case, where the Riemann zeta function is replaced
with a {\itshape Shintani zeta function}. 

These two counting techniques, in combination with the inclusion-exclusion sieve,
are the starting points of the geometric and analytic proofs of Theorems \ref{thm_rc} and \ref{thm_rc_torsion},
respectively!

We note one additional point related to \eqref{eqn_n2}. Wright \cite{W_iwasawa} observed that
this Dirichlet series has the beautiful representation
\begin{equation}\label{eqn_w_iwasawa}
\sum_{[K : \Q] = 2} |\Disc(K)|^{-s} 
=
\prod_p \bigg( \frac{1}{2} \sum_{[K_v : \Q_p] \leq 2} |\Disc(K_v)|_p^s \bigg),
\end{equation}
and proved this in a much more general framework. (He obtained similar formulas
for degree $n$ cyclic extensions of any number field.) His results follow from considering
a nontrivial `twist' of the adelic zeta function of Tate's thesis \cite{tate}. This twist makes
the affine line into a {\itshape prehomogeneous vector space}, with the action $\phi$ of $\GL(1)$
given by $\phi(t) x = t^n x$. 

Essentially, a vector space is prehomogeneous if it has an action of an
algebraic group $G$, which is transitive over $\C$ apart from the vanishing locus of an irreducible polynomial.
This ``prehomogeneous'' property is essential both in Bhargava's work and in the zeta
function approach pioneered by Sato and Shintani \cite{SS}; cubic, quartic, and quintic fields
are parameterized by lattice points\footnote{Not {\itshape all} of the $G(\mathbb{Z})$-orbits 
correspond to fields, as we will see in the cubic case.} up to the action of $G(\mathbb{Z})$,
and these may be counted geometrically or analytically.

Wright's work on \eqref{eqn_w_iwasawa}
and its generalizations mirrors his work with Datskovsky \cite{DW1, DW2, DW3}
on the Shintani zeta function associated to cubic fields, and this
latter work is at the heart of the analytic approach to counting cubic fields. In followup work,
Wright and Yukie \cite{WY} proposed a program to enumerate quartic and quintic fields 
by studying the zeta functions associated to appropriate
prehomogeneous vector spaces. So far their program has not succeeded; however:

{\itshape n = 4, 5.} Bhargava \cite{B_quartic, B_quintic} proved the asymptotic formulas
\begin{equation}\label{eqn_quartic}
N_4(X, S_4) \sim \frac{5}{24} \prod_p \big(1 + p^{-2} - p^{-3} - p^{-4} \big) \cdot X,
\end{equation}
\begin{equation}
N_5(X) \sim \frac{13}{120} \prod_p \big(1 + p^{-2} - p^{-4} - p^{-5} \big) \cdot X.
\end{equation}
The count $N_4(X, S_4)$ includes only quartic fields with Galois group $S_4$, and \eqref{eqn_quartic}
implies
an asymptotic for $N_4(X)$ in combination with work of Baily \cite{baily},
Cohen, Diaz y Diaz, and Olivier \cite{CDO_quartic}, and Wong \cite{wong}. 
Bhargava applies the geometric
approach alluded to in the $N_2(X)$ case (and described later for $N_3(X)$),
and his work may be considered an extension of the methods described in Sections
\ref{sec_dh} and \ref{sec_BST}.

{\itshape n $>$ 5.} For $n > 5$, Bhargava 
conjectured \cite{B_deg_n} that $N_n(X, S_n) \sim C_n X$. The explicit constants $C_n$
have natural interpretations as Euler products; each
Euler factor represents the `probability' that a field $K$ of degree $n$ has a certain
localization. Ellenberg and Venkatesh \cite{EV} have proved the best upper bound to date; namely
$N_n(X) \ll X^{\exp(C \sqrt{\log n})}$.

It is unclear whether prehomogeneous vector spaces may be used to obtain formulas for $N_n(X)$ for $n > 5$.
By a classification theorem of Sato and Kimura \cite{SK}, it is known that there are no reduced, irreducible, and reductive
prehomogeneous vector spaces parameterizing
rings of rank $> 5$; it is not yet known whether there are any prehomogeneous vector spaces not satisfying these assumptions
which parameterize such rings.

\subsection{Counting torsion elements in class groups}
The problem of counting 3-torsion elements in class groups of quadratic fields is related to that of counting
cubic fields. If $D \neq 1$ is a fundamental discriminant, then class field theory provides a bijection between
subgroups of $\Cl(\Q(\sqrt{D}))$ of index 3, and unramified cyclic cubic extensions of $\Q(\sqrt{D})$. Such extensions,
in turn, correspond to (non-Galois) cubic extensions of $\Q$ whose discriminant is $D$. Therefore, counting
3-torsion elements in quadratic fields is equivalent to counting cubic fields whose discriminant is fundamental,
which are those not totally ramified at any prime.

Therefore, this problem may be treated as a variation of the problem of counting cubic fields. This is the approach
of Sections \ref{sec_dh}, \ref{sec_TT}, and \ref{sec_BST}; however, there are other approaches as well, such as that
described in Section \ref{sec_hough}. A general heuristic for statistics of torsion elements in class groups was proposed
by Cohen and Lenstra \cite{CL}; their heuristic has inspired a great deal of interesting followup work, but for the most 
part the Cohen-Lenstra heuristics remain wide open.

\section{Davenport-Heilbronn, Delone-Faddeev, and the main terms}\label{sec_dh}
In this section we discuss the proofs of the main terms in Theorems \ref{thm_rc} and \ref{thm_rc_torsion}.
To summarize, the idea is as follows: We count cubic {\itshape fields} by counting their {\itshape maximal orders},
and so we count all cubic orders and then sieve for maximality. We count cubic orders by counting cubic 
{\itshape rings} and then subtracting the contribution of the reducible rings, and these correspond to binary cubic
{\itshape forms}. (The space of binary cubic forms is {\itshape prehomogeneous} in the sense described previously.) Thus at last we have a geometric lattice-point counting problem, which can be solved by
``brute force'', so to speak. This interpretation in terms of lattice points also opens the door to the use of zeta
functions.

We recommend the paper of Bhargava, Shankar, and Tsimerman \cite{BST} for the simplest self-contained
account of the Davenport-Heilbronn theorems, with complete proofs.

We begin by defining cubic rings and cubic forms.
A {\itshape cubic ring} is a commutative ring which is free of rank 3 as
a $\Z$-module. The discriminant of a cubic ring is defined to be the determinant of the trace form
$\langle x, y \rangle = \Tr(xy)$, and the discriminant of the maximal order of a cubic field is equal to
the discriminant of the field.

The lattice of {\itshape integral binary cubic forms} is defined by
\begin{equation}\label{def_vz}
V_{\Z} := \{ a u^3 + b u^2 v + c u v^2 + d v^3 \ : a, b, c, d \in \Z \},
\end{equation}
and the {\itshape discriminant} of such a form is given by the usual equation
\begin{equation}\label{eqn_disc_formula}
\Disc(f) = b^2 c^2 - 4 a c^3 - 4 b^3 d - 27 a^2 d^2 + 18 abcd.
\end{equation}
There is a natural action of $\GL_2(\Z)$ (and also of $\SL_2(\Z)$) on $V_{\Z}$, given by
\begin{equation}
(\gamma \cdot f)(u, v) = \frac{1}{\det \gamma} f((u, v) \cdot \gamma).
\end{equation}
We call a cubic form $f$ {\itshape irreducible} if $f(u, v)$ is irreducible as a polynomial over $\Q$, and
{\itshape nondegenerate} if $\Disc(f) \neq 0$.

Cubic rings are related to cubic forms by the following correspondence of Delone and Faddeev, 
as further extended by Gan, Gross, and Savin \cite{GGS}:

\begin{theorem}[\cite{DF, GGS}]\label{thm_df} There is a natural, discriminant-preserving 
bijection between the set of $\GL_2(\Z)$-equivalence classes
of integral binary cubic forms and the set of isomorphism classes of cubic rings. Furthermore, under this
correspondence, irreducible cubic forms correspond to orders in cubic fields.

Finally, if $x \in V_{\Z}$ is a cubic form corresponding to a cubic ring $R$, we have
$\Stab_{\GL_2(\Z)}(x) \simeq \Aut(R)$.
\end{theorem}

It is therefore necessary to exclude reducible and nonmaximal rings. The nonmaximality condition is the more 
difficult of the two 
to handle, and for this Davenport and Heilbronn established the following criterion:

\begin{proposition}[\cite{DH, BST}\label{def_up}]
Under the Delone-Faddeev correspondence, a cubic ring $R$ is maximal if any only if its corresponding
cubic form $f$ belongs to the set $U_p \subset V_{\Z}$ for all $p$, defined by the 
following equivalent conditions:
\begin{itemize}
\item
The ring $R$ is not contained
in any other cubic ring with index divisible by $p$.
\item
The cubic form $f$ is not a multiple of $p$, and there is no $\GL_2(\Z)$-transformation of
$f(u, v) = a u^3 + b u^2 v + c u v^2 + d v^3$ such that $a$ is a multiple of $p^2$ and $b$ is a multiple
of $p$.
\end{itemize}
\end{proposition}
In particular, the condition $U_p$ only depends on the coordinates of $f$ modulo $p^2$.

Therefore one can count cubic fields as follows. One obtains an
asymptotic formula for the number of cubic rings of bounded discriminant by counting lattice points
in fundamental domains for the action of $\GL_2(\Z)$, bounded by the constraint $|\Disc(x)| < X$.
The fundamental domains may be chosen so that almost all reducible rings correspond to
forms with $a = 0$, and so these may be excluded from the count.

One then multiplies this asymptotic by the 
product of all the local densities of the sets $U_p$. This gives
a heuristic argument for the main term in \eqref{conj_rc}, and one incorporates a sieve to obtain a proof.

The main term of Theorem \ref{thm_rc_torsion} may be proved similarly. By class field theory,
there is a bijection between subgroups of $\Cl_3(\Q(\sqrt{D}))$ of index 3, and cubic fields of discriminant $D$,
which are precisely those which are not totally ramified at any prime. This ramification condition may
also be detected by reducing cubic forms modulo $p^2$, and thus may be incorporated into Proposition
\ref{def_up}. The remainder of the proof is the same. 

\subsection{The work of Belabas, Bhargava, and Pomerance}\label{subsec_BBP}
In \cite{BBP}, Belabas, Bhargava, and Pomerance (BBP) introduced improvements to Davenport and Heilbronn's
method, and obtained an error term of $O(X^{7/8 + \epsilon})$ in \eqref{conj_rc} and \eqref{eqn_torsion1}. They began by observing that
\begin{equation}\label{eqn_BBP}
N_3^{\pm}(X) = \sum_{q \geq 1} \mu(q) N^{\pm}(q, X),
\end{equation}
where $N^{\pm}(q, X)$ counts the number of cubic orders of discriminant $0 < \pm D < X$ which
are nonmaximal at every prime dividing $q$. This simple observation has proved to be quite useful!
For large $q$, BBP prove that
$N^{\pm}(q, X) \ll X 3^{\omega(q)} / q^2$
using reasonably elementary methods. Therefore, one has
\begin{equation}\label{eqn_BBP2}
N_3^{\pm}(X) = \sum_{q \leq Q} \mu(q) N^{\pm}(q, X) + O(X/Q^{1 - \epsilon}).
\end{equation}

For small $q$, BBP estimate $N^{\pm}(q, X)$ with explicit error terms using geometric methods.
These error terms are good enough to allow them to take the sum in \eqref{eqn_BBP} up to $(X \log X)^{1/8}$, 
which yields a final error term of $O(X^{7/8 + \epsilon})$.

In addition, their methods extend to counting $S_4$-{\itshape quartic} fields, where they obtain a main term of $C_4 X$ 
with error $\ll X^{23/24 + \epsilon}$.

\section{The ``Four Approaches''}\label{sec_four_approaches}
This, then, brings us to Theorems \ref{thm_rc} and \ref{thm_rc_torsion}.
Theorem \ref{thm_rc} was conjectured in print by Roberts \cite{R}, and
implicitly in the foundational work of Datskovsky and Wright \cite{DW1, DW2, DW3}. Roberts gave substantial
and convincing numerical evidence for his conjecture, as well as the heuristic argument described in Section \ref{sec_TT}.

In this paper we will discuss four very different approaches, all independent, to these secondary terms.

{\bf Shintani zeta functions} (Taniguchi-T. \cite{TT, TT_L}).  Taniguchi and the author further
developed the theory of Shintani zeta functions, and proved Theorems \ref{thm_rc} and \ref{thm_rc_torsion}
with error terms of $O(X^{7/9 + \epsilon})$ and $O(X^{18/23 + \epsilon})$, respectively. These proofs closely
follow the heuristic arguments of Roberts and Datskovsky-Wright.

The Shintani zeta function approach is quite flexible, and \cite{TT} contains several generalizations
of the main theorems. Perhaps the most interesting of these is a generalization to counting cubic field discriminants
(or 3-torsion elements of class groups) in arithmetic progressions, where a certain ``non-equidistribution'' phenomenon
arises. For example, when counting cubic field discriminants
$\equiv a \ (\textmod \ 7)$, the secondary term in Theorem \ref{thm_rc} is different for every residue class $a \ (\textmod \ 7)$.

{\bf A refined geometric approach} (Bhargava, Shankar, and Tsimerman \cite{BST}). Bhargava, Shankar,
and Tsimerman gave the first proof of Theorem \ref{thm_rc}, with an error term of $O(X^{13/16 + \epsilon})$. Their proof
follows Davenport and Heilbronn's original work, counting lattice points via geometric arguments, with 
substantial simplifications and improvements. Indeed, these improvements are already evident in their
proofs of the classical Davenport-Heilbronn theorems.

Their proof of the secondary term is based on a ``slicing'' argument which we describe briefly in Section \ref{sec_BST}.
Their treatment of the error terms introduces a certain correspondence for nonmaximal cubic rings;
this correspondence has also proved useful in the context
of Shintani zeta functions, and a combined approach (\cite{BTT}, in progress) has yielded an error term of $O(X^{2/3 + \epsilon})$. 

{\bf Equidistribution of Heegner points} (Hough \cite{H}). Hough obtained a statement closely related to
Theorem \ref{thm_rc_torsion}. His scope is limited to counting 3-torsion elements in imaginary quadratic fields $\Q(\sqrt{D})$
with $D \equiv 2 \ (\textmod \ 4)$; in this context, he proves a version of Theorem \ref{thm_rc_torsion}, but with an error term
larger than $X^{5/6}$. In followup work (in preparation), he proves the secondary term for a {\itshape smoothed} variation
of this sum.

Hough derives this result as a consequence of an equidistribution theorem for the Heegner points
associated to the 3-part of class groups of imaginary quadratic fields; his proof of the Davenport-Heilbronn
theorem appears as a bonus, and without reference to binary cubic forms. His techniques also apply to the $k$-parts of these class groups, 
for odd $k \geq 5$, where he obtains secondary terms (of order $X^{1/2 + 1/k}$). His error terms
are currently larger than both the secondary and the main terms for $k \geq 5$, but this work sheds light on a notoriously
difficult problem.

{\bf Trigonal curves and the Maroni invariant} (Zhao \cite{Z}). Zhao enumerates cubic extensions of the rational function field
$\Fqt$ using algebraic geometry. Such extensions are in bijection with isomorphism classes of smooth trigonal curves, that is, smooth 3-fold covers of
$\PP^1$, and these may be counted by embedding them in certain surfaces $F_k$. He obtains a secondary term as a consequence
of a bound for the integer $k$, called the Maroni invariant; for now his error terms are larger than this secondary term.
\\
\\
In all of these approaches, the secondary term is easier to see than it is to prove. In particular, each of these approaches
yield the secondary term in a natural way, but it is not {\itshape a priori} evident that the error terms can be made smaller
than $X^{5/6}$!

In what follows, we will describe, insofar as we can in a couple of pages, why each of these approaches naturally yields a secondary term.
For the Shintani zeta function approach, we will say only a little about the error term, and for the
other three approaches we will say nothing. Indeed, we will use the notation $O(\cdots)$ whenever 
the details
are best left to the respective papers.

\section{The Shintani Zeta Function Approach}\label{sec_TT}

Taniguchi and the author \cite{TT} proved Theorems \ref{thm_rc} and \ref{thm_rc_torsion} using the analytic theory of Shintani
zeta functions. The {\itshape Shintani zeta functions} associated to the space of binary cubic forms
are defined\footnote{The Shintani zeta functions are commonly defined in terms of $\SL_2(\Z)$
rather than $\GL_2(\Z)$, multiplying \eqref{eqn_shintani_def_1}
by a factor of 2.} 
by the Dirichlet series
\begin{equation}\label{eqn_shintani_def_1}
\xi^{\pm}(s) := \sum_{x \in \GL_2(\Z) \backslash V^{\pm}_{\Z}} \frac{1}{|\Stab_{\GL_2(\Z)}(x)|} |\Disc(x)|^{-s}
= \sum_{\substack{\pm \Disc(R) > 0}} \frac{1}{|\Aut(R)|} |\Disc(R)|^{-s}.
\end{equation}
In the former sum, $V_{\Z}$ is the lattice
defined in \eqref{def_vz}, the sum is over elements of positive or negative discriminant respectively,
and $\Stab(x)$ is the stabilizer
of $x$ in $\GL_2(\Z)$. The latter sum is over isomorphism classes of cubic rings.
The equality of these two sums follows from the Delone-Faddeev
correspondence (Theorem \ref{thm_df}).

This is an example of a zeta function associated to a prehomogeneous vector space; in this case
the word `prehomogeneous' reflects the fact that $\GL_2(\R)$ acts transitively on the positive- and
negative-discriminant loci $V^{\pm}_{\R}$. Sato and Shintani \cite{SS} developed a general theory
of zeta functions associated to prehomogeneous vector spaces, and for the space of binary cubic forms
Shintani proved \cite{S} 
that these zeta functions enjoy analytic continuation and an explicit functional equation.
These zeta functions have poles at $s = 1$, and, anomalously, at $s = 5/6$.

Shintani's work opens the door to the study of cubic fields using analytic number theory. In particular, by Perron's formula and 
standard techniques
we have
\begin{equation}\label{eqn_perron1}
\sum_{\substack{x \in \GL_2(\Z) \backslash V_{\Z} \\ \pm \Disc(x) < X}} \frac{1}{|\Stab(x)|}
= \int_{2 - i \infty}^{2 + i \infty} \xi^{\pm}(s) \frac{X^s}{s} ds
= \Res_{s = 1} \xi^{\pm}(s) X + \frac{6}{5} \Res_{s = 5/6} \xi^{\pm}(s) X^{5/6} + O(X^{3/5 + \epsilon}).
\end{equation}
Although the left side is not the counting function of cubic fields, for the first time {\itshape we see the $X^{5/6}$ secondary term.}
This, therefore, gives an explanation for the secondary terms in Theorems \ref{thm_rc} and \ref{thm_rc_torsion}:
because the Shintani zeta functions have secondary poles.\footnote{
Of course, this begs the question of why the Shintani zeta functions have poles at $s = 5/6$. This can of course be explained
by the various calculations in Shintani's work, but for now a satisfying highbrow argument is unknown to the author.
}

To count cubic fields, the most important step is sieving for maximality. This is possible by work of Datskovsky and
Wright \cite{DW1, DW2, DW3}, who gave Shintani's work an adelic formulation along the lines of Tate's thesis \cite{tate}.
This allowed them to incorporate a variety of conditions into the definition of the Shintani zeta functions, including
the Davenport-Heilbronn maximality conditions (given in Proposition \ref{def_up}). 

A bogus proof of Roberts' conjecture is as follows. For a set of primes $\mathcal{P}$,
define the {\itshape $\mathcal{P}$-maximal Shintani zeta function} by the Dirichlet series \eqref{eqn_shintani_def_1},
with the added condition that the only cubic forms counted are those in the set $U_p$ for each $p \in \calP$. Therefore,
letting $\mathcal{P}$ be the set of primes $< X$, it follows that 
\begin{equation}\label{eqn_rc_heur}
\sump_{\pm \Disc(x) < R} \frac{1}{|\Aut(R)|}
= \int_{2 - i \infty}^{2 + i \infty} \xi_{\mathcal{P}}^{\pm}(s) \frac{X^s}{s} ds
= \Res_{s = 1} \xi_{\mathcal{P}}^{\pm}(s) X + \frac{6}{5} \Res_{s = 5/6} \xi_{\mathcal{P}}^{\pm}(s) X^{5/6} + O(X^{3/5 + \epsilon}),
\end{equation}
where the sum is over all {\itshape maximal} cubic rings\footnote{Reducible maximal cubic rings correspond to fields
of degree $\leq 2$, which can be separately counted and therefore subtracted with good error terms.}, or equivalently cubic fields,
with $\pm \Disc(R) < X$. Unfortunately 
the error term depends on the particular zeta function $\xi_{\mathcal{P}}$, and therefore the implied constant depends on $X$
(indeed, when standard techniques are used, exponentially so!)

However, this flawed argument enabled Datskovsky and Wright \cite{DW3} to give a rigorous analytic proof of
the Davenport-Heilbronn theorem. The key step\footnote{We have adapted their argument somewhat to our point of view.
Datskovsky and Wright also need, as did Davenport and Heilbronn, the bound \eqref{eqn_BBP2} for the number of cubic rings 
nonmaximal at any prime $> Y$.}
is to use \eqref{eqn_rc_heur}, but with $\mathcal{P}$ equal to the set of all primes $< Y$, where $Y$ tends to infinity
very slowly with $X$. Given their adelic framework, their work yielded
an asymptotic formula for the number of cubic extensions of any global
field (of characteristic not 2 or 3). 

Moreover, Roberts \cite{R} computed the limit of the secondary terms in \eqref{eqn_rc_heur} and found an excellent
match with numerical data. This led him to conjecture that \eqref{eqn_rc_heur} (equivalently, \eqref{conj_rc}) was correct 
with some undetermined error term smaller than  $X^{5/6}$.

He closed his paper with the following (paraphrased very slightly):

\begin{quote}
The pessimistic discussion in \cite{DW3} suggests to us that the way may be difficult. However,
one ingredient of a proof might be the functional equation of $\xi^{\pm}_{\calP}(s)$
with respect to $s \rightarrow 1 - s$, studied in \cite{DW1, DW2}. Another
ingredient might be [\cite{SS}, Thm. 3], which concerns
growth of arithmetic functions whose associated Dirichlet series satisfy
such a functional equation.
\end{quote}

This is actually a good summary of our proofs of Theorems \ref{thm_rc} and \ref{thm_rc_torsion}! To get the ball rolling, we needed to incorporate
the sieve used by Belabas, Bhargava, and Pomerance described in Section \ref{subsec_BBP}. In place of the $\calP$-maximal
zeta function we introduced the {\itshape $q$-nonmaximal zeta function}, which counts only cubic rings {\itshape not} maximal at each prime dividing
$q$. This allowed us to count maximal cubic rings, and the irreducible such rings correspond to cubic fields.

\begin{remark}
We count 3-torsion elements of class groups by a variation of this argument. Such elements correspond to cubic 
fields not totally ramified at any prime, and so we expand the definition of the $q$-nonmaximal zeta function to count cubic rings
either nonmaximal or totally ramified at every prime dividing $q$.
\end{remark}

The equation \eqref{eqn_rc_heur} holds for each $q$-nonmaximal zeta function, 
and the $q$-dependence of the error term can be bounded in terms of the exponential sum
\begin{equation}\label{def_phihat}
\sum_{x \in V_{\Z / q^2 \Z}} |\widehat{\Phi}_q(x)| := \sum_{x \in V_{\Z / q^2 \Z}} \bigg|
\frac{1}{q^8} \sum_{y \in V_{\Z / q^2 \Z}} \Phi_q(y) \exp(2 \pi i [x, y] / q^2) \bigg|,
\end{equation}
where $\Phi_q(x)$ is the characteristic function of the Davenport-Heilbronn $q$-nonmaximality condition, and
$[x, y] = x_4 y_1 - \frac{1}{3} x_3 y_2 + \frac{1}{3} x_2 y_3 - x_1 y_4.
$
Each inner sum appears naturally when we shift the contour in \eqref{eqn_rc_heur} to the left of $\Re(s) = 0$ and apply
the functional equation. Trivially, the sum in \eqref{def_phihat} is bounded by $q^8$; this is enough to allow us to take $Q = X^{1/25}$ in \eqref{eqn_BBP2}
and obtain an error term of $X^{24/25 + \epsilon}$ in the Davenport-Heilbronn theorems.

In \cite{TT} we incorporated several improvements to lower the error terms in
Theorems \ref{thm_rc} and \ref{thm_rc_torsion} to $O(X^{7/9 + \epsilon})$ and $O(X^{18/23 + \epsilon})$, respectively.
Among these is a careful analysis of the sum in \eqref{def_phihat}, carried out in \cite{TT_L}. For each
$x$ we obtained exact formulas for the inner sum in \eqref{def_phihat}, proving that the outer sum is $\ll q^{1 + \epsilon}$.

\subsection{Non-equidistribution in arithmetic progressions}\label{sec_TT_ap}
The Shintani zeta function approach can also be used to study the distribution of cubic field
discriminants in arithmetic progressions. To do this, one twists the Shintani zeta functions
by Dirichlet characters, obtaining $L$-functions which are also proved to 
enjoy analytic continuation and a functional
equation.

The details are complicated, so here we simply present representative
numerical data. The tables below list the number of cubic fields $K$ with
$0 < \Disc(K) < 2 \cdot 10^6$ and $\Disc(K) \equiv a \ (\textmod \ m)$, for $m = 5$ and $m = 7$, and each
residue class $a$.

\begin{center}
\begin{tabular}{l | c | c | c | c | c  }
Discriminant modulo 5 & 0 & 1 & 2 & 3 & 4  \\ \hline
Actual count & 21277 & 22887 & 22751 & 22748 & 22781 \\
Theoretical result & 21307 & 22757 & 22757 & 22757 & 22757 \\
Difference & 30 & 130 & 6 & 9 & 24 \\
\end{tabular}
\end{center}

\begin{center}
\begin{tabular}{l | c | c | c | c | c | c | c}
Discriminant modulo 7 & 0 & 1 & 2 & 3 & 4 & 5 & 6 \\ \hline
Actual count & 15330 & 17229 & 14327 & 15323 & 17027 & 18058 & 15150 \\
Theoretical result & 15316 & 17209 & 14277 & 15316 & 17024 &  18063 & 15131 \\
Difference & 14 & 20 & 50 & 7 & 3 & 5 & 19 \\
\end{tabular}
\end{center}
The data modulo 5 is essentially equidistributed, apart from the $a = 0$ column. 
(The relative deficit of cubic fields divisible by 5, or by any other modulus, 
can be explained by a careful reading of Davenport and Heilbronn's
original paper.) However, the data modulo 7 shows a striking lack of equidistribution.\footnote{
The counting functions for cubic field discriminants $\equiv 0 \ (\textmod \ 7)$ and $\equiv 3 \ (\textmod \ 7)$
are not the same; the equality of the theoretical results above is a coincidence of roundoff error.}

The theoretical results above come from a generalization of Theorem 
\ref{thm_rc} to arithmetic progressions; the secondary term now involves a
sum of residues of twisted Shintani zeta functions. These residues are evaluated in 
\cite{TT_L}, and for characters $\chi$ for which $\chi^6 \neq 1$, the associated residue
at $s = 5/6$ is zero. However, when $\chi^6 = 1$ this residue is often nonzero,
and nonzero residues explain biases in arithmetic progressions such as the one
seen above.

\section{A Refined Geometric Approach}\label{sec_BST}

To give an overview of Bhargava, Shankar, and Tsimerman's proof \cite{BST}, we start with
their treatment of the main term in Theorem \ref{thm_rc}.
Later, we will describe how the secondary term appears in their work, but only
in the context of irreducible cubic forms (corresponding to cubic orders). 

The argument in \cite{BST} largely follows Davenport and Heilbronn's original work.
Fix a sign, and let $\Virr$ denote the irreducible points $x \in V_{\Z}$ with $\pm \Disc(x) > 0$.
Write $n_+ = 6$, $n_- = 2$ for the order of the stabilizers of the action of $\GL_2(\R)$ on $V_{\R}^{\pm}$.
The equality $n_+ = 3 n_-$ reflects the fact that cubic fields of negative discriminant are three times as common 
as those of positive discriminant.

Write $\calF$ for a certain fundamental domain for $\GL_2(\Z) \backslash \GL_2(\R)$
in $\GL_2(\R)$ (see (12) of \cite{BST}), originally constructed by Gauss. For any 
$v \in V^{\pm}$, consider the multiset $\calF v$, where the multiplicity of an element $x \in V^{\pm}$ is equal
to the number of $g \in \calF$ for which $g v = x$. It is readily checked that 
each element $x \in G_{\Z} \backslash V^{\pm}_{\Z}$ is represented 
in this multiset $n_{\pm}/m(x)$
times, where $m(x)$ denotes the size of the stabilizer of $x$ in $\GL_2(\Z)$. We conclude that 
\begin{equation}\label{eqn_BST1}
N(V^{\pm}; X) := \sum_{\substack{0 < \pm \Disc(x) < X \\ x \ {\textnormal irred.}}} \frac{1}{|\Stab(x)|} = 
\frac{1}{n_{\pm}} \{ x \in \calF v \cap \Virr \ : \ |\Disc(x)| < X \}.
\end{equation}

The reducible points can be handled without too much difficulty: The number of reducible integral binary
cubic forms $a u^3 + b u^2 v + c u v^2 + d v^3$ in the multiset $\calR_X(v) := \{ w \in \calF v \ : \ |\Disc(w)| < X \}$,
satisfying $a \neq 0$, is $\ll X^{3/4 + \epsilon}$. Conversely, all cubic forms with $a = 0$ are reducible.
This condition $a \neq 0$ meshes well with the geometric arguments that follow.
In addition, the weight $\frac{1}{|\Stab(x)|}$ occurring in \eqref{eqn_BST1} can be neglected, 
as there are $\ll X^{1/2}$ points $x$ with $|\Disc(x)| < X$ with nontrivial stabilizer.

The formula in \eqref{eqn_BST1} does not depend on $v$, and one innovation of \cite{BST} (previously used by Bhargava in \cite{B_quartic}) is to average over many $v$.
In particular, define $B := \{ w = (a, b, c, d) \in V : \ 3a^2  + b^2 + c^2 + 3d^2 \leq 10, \ |\Disc(w)| \geq 1 \}$, and we have
\begin{equation}\label{eqn_BST_avg}
N(V^{\pm}; X) = 
\frac{
\int_{v \in B \cap V^{\pm}} 
\frac{1}{n_{\pm}} \{ x \in \calF v \cap \Virr \ : \ |\Disc(x)| < X \} |\Disc(v)|^{-1} dv
}
{
\int_{v \in B \cap V^{\pm}} 
|\Disc(v)|^{-1} dv
}.
\end{equation}
The authors describe this step as ``thickening the cusp''. As $|\Disc(v)|^{-1} dv$ is $\GL_2(\R)$-invariant, we may rewrite this as
\begin{equation}
N(V^{\pm}; X) = \frac{1}{M^{\pm}} \int_{g \in \calF}
\# \{ x \in \Virr \cap gB : |\Disc(x)| < X \} dg,
\end{equation}
where $M^{\pm} := \frac{n_{\pm}}{2\pi} \int_{v \in B \cap V^{\pm}} |\Disc(v)|^{-1} dv.$ This is rewritten as
\begin{equation}\label{eqn_bst_main}
N(V^{\pm}, X) = \frac{1}{M^{\pm}}
\int_{g \in N'(a) A' \Lambda} \# \{x \in \Virr \cap B(n, t, \lambda, X) \} t^{-2} dn d^{\times} t d^{\times} \lambda
\end{equation}
using a standard decomposition of $\calF$ ((12) of \cite{BST}), where $B(n , t, \lambda, X)$
is the region $\{ x \in g B:  |\Disc(x)| < X \}$. A proposition of Davenport establishes that in this
situation, the count of lattice points above is well approximated by the volume of $B(n, t, \lambda, X)$,
provided that $t$ is not too large. Conversely, when $t$ is large all the lattice points are in the cusp $a = 0$,
and hence reducible and not counted. Thus, 
the above is reduced to 
\begin{equation}\label{eqn_bst_main_red}
N(V^{\pm}, X) =\frac{1}{M^{\pm}}
\int_{\substack{g \in N'(a) A' \Lambda \\ t < C^{1/3} \lambda^{1/3}}} \Vol(B(n, t, \lambda, X))t^{-2} dn d^{\times} t d^{\times} \lambda + O(\cdots).
\end{equation}

The condition on $t$ is now removed, subject to an error term, and this integral is equal to
\begin{equation}
N(V^{\pm}, X) = \frac{1}{2 \pi M^{\pm}} \int_{v \in B \cap V^{\pm}} \Vol(\calR_X(v)) |\Disc(v)|^{-1} dv + O(\cdots).
\end{equation}
This volume does not in fact depend on $v$, 
so that we have (excluding error terms)
\begin{equation}\label{eqn_bst_main_term}
N(V^{\pm}, X) = \frac{1}{n_{\pm}} \Vol(\calR_X(v)) = \frac{\pi^2}{12 n_{\pm}} X.
\end{equation}

\subsection{Origin of the secondary term}
We now explain how Bhargava, Shankar, and Tsimerman refine these calculations to obtain a secondary term.
Our brief explanation will necessarily be somewhat vague, which should encourage the reader
to read \cite{BST}. We focus on their count of irreducible binary cubic {\itshape forms}; as in \cite{BBP, TT}, they also incorporate a sieve to
count cubic fields.

They begin by incorporating several tweaks to \eqref{eqn_bst_main}. For technical reasons, they count discriminants in dyadic intervals
$[X/2, X]$. In addition, they introduce a smooth function $\Psi_0(t)$ on $\R_{\geq 0}$, such that $\Psi_0(t) = 0$ for $t \leq 2$ and $\Psi_0(t) = 1$
for $t \geq 3$. They thus rewrite \eqref{eqn_bst_main} (restricted to a dyadic interval $[X/2, X]$) as
\begin{multline}\label{eqn_bst_main_2}
N(V^{\pm}, [X/2, X]) = \frac{1}{M^{\pm}}
\int_{g \in N'(a) A' \Lambda} \bigg(\Psi_0\bigg( \frac{t \kappa}{\lambda^{1/3}} \bigg) + \Psi\bigg( \frac{t \kappa}{\lambda^{1/3}} \bigg) \bigg)
\times \\ \# \{x \in \Virr \cap B(n, t, \lambda, [X/2, X]) t^{-2} dn d^{\times} t d^{\times} \lambda,
\end{multline}
where $\Psi := 1 - \Psi_0$, and $\kappa$ is a parameter to be chosen later (to minimize error terms). The decomposition $\Psi + \Psi_0 = 1$ splits this integral into two, and we restrict our attention to the 
$\Psi_0$ part, which yields the secondary term.

They now ``slice'' the count of binary cubic forms by the first coordinate. In particular, for $a \in \mathbb{Z}$,
let $B_a(n, t, \lambda, [X/2, X])$ denote the set of binary cubic forms in $B_a(n, t, \lambda, [X/2, X])$
whose $u^3$ coefficient is equal to $a$. Then, we have
\begin{equation}
\# \big\{ x \in \Virr \cap B(n, t, \lambda, [X/2, X]) \big\} = \sum_{\substack{a \in \Z \\ a \neq 0}}
\# \big\{ x \in \Virr \cap B_a(n, t, \lambda, [X/2, X]) \big\}.
\end{equation}
Skipping ahead in their work, the $\Psi_0$-contribution to \eqref{eqn_bst_main_2} is equal to
\begin{equation}\label{eqn_BST_pre_mellin}
\frac{2}{3 M^{\pm}} \sum_{a = 1}^{\infty} \int_{\lambda = (\sqrt{3 / 2})^3 / C}^{X^{1/4}} \int_{u > 0} \Phi_0\bigg(\frac{u^{1/3} \kappa}{a^{1/3}} \bigg)
\frac{\lambda^{10/3} u^{1/3}}{a^{1/3}} \Vol\big(B_u([X/(2 \lambda^4), X/\lambda^4])\big) d^{\times} u d^{\times} \lambda,
\end{equation}
where $B_a([Y/2, Y])$ denotes the set of cubic forms in $B$ with first coordinate $a$ and discriminant in $[Y/2, Y]$. By Mellin inversion,
we have
\begin{equation}\label{eqn_BST_mellin}
\sum_{a = 1}^{\infty} a^{-1/3} \Phi_0\bigg( \frac{u^{1/3} \kappa}{a^{1/3}} \bigg) = 
\zeta(1/3) + 3 \widetilde{\Phi}_0(-2) (\kappa^3 u)^{2/3} + O(\cdots),
\end{equation}
and it is this (negative) $\zeta(1/3)$ which contributes the secondary term. In contrast, the second term of \eqref{eqn_BST_mellin}
(which ``looks larger'', especially in light of the eventual choice $\kappa = X^{1/12}$) 
is
reincorporated into \eqref{eqn_BST_pre_mellin} and then combined with the $\Psi$ term of \eqref{eqn_bst_main_2} to obtain the main term
in \eqref{eqn_bst_main_term}. When the $\zeta(1/3)$ term is plugged into \eqref{eqn_BST_pre_mellin}, the resulting integral
can be evaluated without undue difficulty, and it has order $X^{5/6}$.

\begin{remark} As Shankar has explained to the author, a more geometric argument can also be given, which does not 
rely on the analytic continuation of the zeta function. However, it is unclear whether the resulting error term 
would be smaller than $X^{5/6}$, when combined with the sieve.
\end{remark}

\subsection{A correspondence for cubic forms}\label{sec_corr}
To count cubic fields one must sieve for maximality, using
the Davenport-Heilbronn conditions
of Proposition \ref{def_up}. These could be treated in a naive manner, but the authors improve their
error terms by introducing a useful correspondence for nonmaximal cubic forms.\footnote{As Shankar explained to me, in the case of cubic fields it is this
correspondence which allows them to cross the $X^{5/6}$ barrier.}

We introduce some notation. Let $N^{\pm}(V_{\Z}; X)$ count all cubic forms $v$ with $0 < \pm \Disc(v) < X$.
For a prime $p$, let $N^{\pm}(U_p^c; X)$ count those cubic forms which are nonmaximal at $p$, i.e., which do not lie
in the set $U_p$ defined
in Proposition \ref{def_up}. (In referring to cubic forms as ``nonmaximal'' we are implicitly appealing to the Delone-Faddeev correspondence.)
Finally, for any $\alpha \in \PP^1(\Z/p\Z)$, let $N^{\pm}(V_{p, \alpha}; X)$ count those cubic forms $v$ such that the reduction of $v$ modulo
$p$ has $\alpha$ as a root. 

With this notation, we have the following correspondence:
\begin{proposition}\cite{BST} We have the identity
\begin{equation}\label{eqn_bst_id}
N^{\pm}(U_p^c; X) = \sum_{\alpha \in \PP^1(\F_p)} N^{\pm}(V_{p, \alpha}; X/p^2)
- \sum_{\alpha \in \PP^1(\F_p)} N^{\pm}(V_{p, \alpha}; X/p^4) + N^{\pm}(V_{\Z}; X/p^4).
\end{equation}
\end{proposition}
This identity generalizes in a straightforward way from prime $p$ to squarefree $q$, where
we count cubic forms not in $U_p$ for any $p | q$. The identity amounts to a combinatorial argument, enumerating ways 
in which a nonmaximal cubic ring is contained in a larger ring (with smaller discriminant), and then counting these larger rings.

This identity is crucial in \cite{BST}, and it can also be translated into an identity for the $q$-nonmaximal Shintani zeta function!
In work in progress by Bhargava, Taniguchi, and the author \cite{BTT} we are developing this combined approach. We have
proved an error term of $O(X^{2/3 + \epsilon})$ in Theorem \ref{thm_rc} as well as a secondary term for relative cubic extensions
of quadratic base fields, and we are currently working on further extensions and generalizations.

\section{Equidistribution of Heegner Points}\label{sec_hough}

Remarkably, Hough \cite{H} demonstrated that the Delone-Faddeev correspondence is not the only path to the Davenport-Heilbronn
theorem. Hough studied the distribution of Heegner points associated to 3-torsion elements of ideal class groups of quadratic
fields, and he 
obtained a result related to Theorem \ref{thm_rc_torsion} as a consequence.

For simplicity, Hough worked only with imaginary discriminants $D \equiv 2 \ (\textmod \ 4)$; note that any such discriminant
is not fundamental, and the discriminant of the field $\Q(\sqrt{D})$ is $4D$.\footnote{It seems that the restrictions that $D \equiv 2 \ (\textmod \ 4)$
and that $D$ be negative could both be lifted with some effort. The restriction on the sign is the more serious of the two, as elements
of class groups of $\Q(\sqrt{D})$ for positive $D$ correspond to closed geodesics rather than Heegner points. However, Duke's result
holds for either sign, so in \cite{H} Hough expresses optimism that the positive discriminant case of Theorem \ref{thm_rc_torsion} could
be handled as well.} He proves that
\begin{equation}
\sum_{\substack{0 < D < X \\ D \equiv 2 \ (\textmod \ 4)}} \# \Cl_3(-D) = \frac{4}{\pi^2} X + O(X^{19/20 + \epsilon}),
\end{equation}
and for ``good''\footnote{Compact support and infinitely differentiable.} test functions $\phi$, he proves that
\begin{equation}\label{eqn_hough2}
\sum_{\substack{0 < D \\ D \equiv 2 \ (\textmod \ 4)}} \# \Cl_3(-D) \phi(D/X)= \frac{4}{\pi^2} \widehat{\phi}(1) X + O(X^{7/8 + \epsilon}).
\end{equation}
The correct secondary term of order $X^{5/6}$ appears in both formulas, in spite of the larger error terms, and in 
in followup work (in progress) he obtains an error less than $X^{5/6}$ 
in \eqref{eqn_hough2}, thereby obtaining
a proof of a related secondary term.

Hough's work has very different prospects for generalization from the approaches described previously. In particular, his
approach can be used to study $k$-torsion in class groups for odd $k > 3$, and he makes the following conjecture:

\begin{conjecture}[Hough \cite{H}]
For good test functions $\phi$ and odd $k \geq 3$, we have
\begin{equation}\label{eqn_hough3}
\sum_{\substack{0 < D \\ D \equiv 2 \ (\textmod \ 4)}} \# \Cl_{k}(-D) \phi(D/X)= \frac{4}{\pi^2} \widehat{\phi}(1) X + 
C_{1, k} \widehat{\phi}\bigg(\frac{1}{2} + \frac{1}{k} \bigg) X^{1/2 + 1/k} + o(X^{1/2 + 1/k}),
\end{equation}
where
\begin{multline}
C_{1, k} := \frac{1}{6 k} \frac{\zeta(1 - \frac{2}{k})}{\zeta(2)} \frac{\Gamma(\frac{1}{2}) \Gamma(\frac{1}{2} - \frac{1}{k})}{\Gamma(1 - \frac{1}{k})}
\Big(1 - 2^{1/k} + 2^{1 - 1/k}\Big) \times \\ 
\prod_{p > 2}
\bigg(1 + \frac{1}{p + 1} \Big( p^{-1/k} - p^{-1 + 2/k} - p^{-1 + 1/k} - p^{-1} \Big) \bigg).
\end{multline}
\end{conjecture}

Indeed he ``proves'' his conjecture, but with error terms larger than $X$.
Even a proof of the main term for the single case $k = 5$ would be a major achievement. Indeed, the
average size of $\Cl_p(-D)$ was conjectured by Cohen and Lenstra \cite{CL} to be 2 for each prime $p \geq 3$, and 
Cohen and Lenstra further conjectured that the $p$-part of $\Cl_p(-D)$ is isomorphic to a fixed $p$-group $H$ with probability proportional
to $\frac{1}{|\Aut(H)|}$. No case of these conjectures is currently known for any $p \geq 5$, and essentially nothing is known for $p = 3$
beyond Theorem \ref{thm_rc_torsion}.

\subsection{Heegner points and equidistribution}
Hough's main result concerns the equidistribution of Heegner points associated to the 3-part of the class group.
We recall the relevant background.
If $\mfa  = [x_1, x_2]$ is an ideal of an imaginary quadratic field $\Q(\sqrt{-D})$,
where by reordering if necessary we have $x_1/x_2 \in \HH := \{ z \in \C \ : \ \Im(z) > 0 \}$, we call
$x_1/x_2$ the {\itshape Heegner point $z_{\mfa}$ associated to the ideal} $\mfa$. A change of basis for $\mfa$
corresponds to a linear fractional transformation $z_{\mfa} \rightarrow \frac{a z_{\mfa} + b}{c z_{\mfa} + d}$ with
$ad - bc = 1$, and therefore the Heegner point $z_{\mfa}$ is well-defined in the quotient $\SL_2(\Z) \backslash \HH$. One
checks that $z_{\mfa} = z_{\alpha \mfa}$ for any $\alpha \in \Q(\sqrt{-D})$, and hence $z_{\mfa}$ depends only on
the class of $\mfa$ in $\Cl(\sqrt{-D})$. One similarly checks that $z_{\mfa} \neq z_{\mfb}$ if $\mfa$ and $\mfb$ are inequivalent
in $\Cl(\sqrt{-D})$.

Therefore, there are $h(-D)$ Heegner points associated to $\Q(\sqrt{-D})$ in $\SL_2(\Z) \backslash \HH$,
and in \cite{duke} Duke proves that they equidistribute with respect to the hyperbolic measure $\frac{dx \ dy}{y^2}$
as $D \rightarrow \infty$. Hough proves the same for the Heegner points associated to the 3-part of the ideal class group:
\begin{theorem}[Hough \cite{H}]
Let $B$ be a bounded Borel measurable subset of $\mathbb{H}$ having boundary of measure zero. Then as $D \rightarrow
\infty$,
\begin{equation}
\sum_{\substack{0 > -D > -X \\ d \equiv 2 \ (\textmod \ 4) \\ \textnormal{squarefree}}}
\#\big\{ \mfa \ \textnormal{primitive, nonprincipal} : [\mfa] \in Cl_3(\Q(\sqrt{-D})), z_{\mfa} \in B\big\}
\sim
\frac{6X}{\pi^3} \vol(B).
\end{equation}
\end{theorem}
Here we take $B \subseteq \mathbb{H}$ rather than $B \subseteq \SL_2(\Z) \backslash \HH$, so that each Heegner point is counted
once for each representative of its $\SL_2(\Z)$ orbit in $B$. (We recall that a {\itshape primitive} ideal is one without any
rational integer divisors other than 1.)

Straightaway one can see that this is an equidistribution statement. But it also implies the
Davenport-Heilbronn theorem! Taking $B$ to be the standard 
(or any other) fundamental domain 
for the action of $\SL_2(\Z)$ on $\HH$,
one counts each Heegner point (and thus
each nontrivial element of $\Cl_3(\Q(\sqrt{-D}))$) exactly once. The Davenport-Heilbronn theorem
follows from the classical computation $\vol(B) = \frac{3}{\pi}$.

We will follow Hough's approach sufficiently far as to see the secondary term. His work begins with the following
parameterization of ideals (and therefore Heegner points), essentially due to Soundararajan \cite{sound}.
\begin{proposition}\cite{H, sound}
Let $D \equiv 2 \ (\textmod \ 4)$ be squarefree and let $k \geq 3$ be odd. The set
\begin{equation}\label{eqn_hough}
\{ (l, m, n, t) \in (\Z^+)^4 \ : \ l m^k = l^2 n^2 + t^2 D, \ l | D, \ (m, ntD) = 1\}
\end{equation}
is in bijection with primitive ideal pairs $(\mfa, \overline{\mfa})$ with $\mfa \neq 1$ and $\mfa^k$ principal
in $\Q(\sqrt{-D})$. Explicitly, the ideal $\mfa$ is given as a $\Z$-module by
$$
\mfa = [lm, lnt^{-1} + \sqrt{-D}]
$$
where $\N \mfa = lm$ and $t^{-1}$ is the inverse of $t$ modulo $m$.
\end{proposition}
The proof is relatively straightforward. For example, if $\mfa^k$ is principal and coprime to $D$, then
$\mfa^k = (n + t\sqrt{-D})$ and taking norms yields $(\N \mfa)^k = n^2 + t^2 D$, a solution to \eqref{eqn_hough}.
Ideals $\mfa$ not coprime to $D$ yield solutions to \eqref{eqn_hough} with $l > 1$. 

As one of the two main steps in his proof, Hough now applies this parameterization
to count Heegner points in the region
\begin{equation}
R_Y := \Big\{ z \in \HH \ : -\frac{1}{2} <  \Re(z) < \frac{1}{2}, \ \ \Im(z) > \frac{1}{Y} \Big\}.
\end{equation}
(This establishes the vertical distribution of Heegner points, and he separately establishes
that they are equidistributed horizontally.)
He proves that the number of Heegner points in $R_Y$ is equal to
\begin{equation}\label{eqn_heeg_pts}
\frac{6}{\pi^3} Y X + C_{5/6} X^{5/6} + O(\cdots),
\end{equation}
for a fairly complicated error term which depends on both $X$ and $Y$. The secondary term appears here for the following reason:
By the parameterization above,
we may write an ideal $\mfa$ as $ [lm, l n t^{-1} + \sqrt{-D}]$. We have $l m^3 = l^2 n^2 + t^2 D$, so that $lm \geq D^{1/3}$,
implying that the Heegner point $z_{\mfa}$ has imaginary part $\frac{\sqrt{D}}{lm} \leq D^{1/6} \leq X^{1/6}$. 
In other words, the Heegner
points equidistribute as $X \rightarrow \infty$, but for fixed $X$ they do not appear high in the cusp.

The hyperbolic volume of the subset of $R_Y$ with $\Im(z) > X^{1/6}$ is equal to 
\begin{equation}
\int_{x = -1/2}^{1/2} \int_{y = X^{1/6}}^{\infty} \frac{dx \ dy}{y^2} = X^{-1/6}.
\end{equation}
Accounting (much more carefully than done here!) 
for the fact that this subset contains no Heegner points yields the secondary
negative term in \eqref{eqn_heeg_pts}, and an analogous argument explains the second term in \eqref{eqn_hough3}.

\begin{remark} Hough's methods also extend to arithmetic progressions, where he finds the same
lack of equidistribution observed in Section \ref{sec_TT_ap}; some results along these lines are in preparation.
\end{remark}

\section{Hirzebruch Surfaces and the Maroni Invariant}\label{sec_zhao}
Zhao's work \cite{Z} is still in preparation, so we offer only a very brief overview here. Using algebraic geometry,
Zhao estimates the number of cubic extensions of the {\itshape rational function field}
$\Fqt$.

Fix a finite field $\mathbb{F}_q$ of characteristic not equal to 2 or 3, and let $S(2N)$ be the number
of isomorphism classes of cubic extensions $K/\Fqt$ of degree 3, such that $|\Disc(F)| = q^{2N}$. In \cite{DW3},
Datskovsky and Wright proved that
\begin{equation}\label{eqn_DW_ff}
S(2N) = 2 \frac{\Res_{s = 1} \zeta_{\Fqt}(1)}{\zeta_{\Fqt}(3)} q^{2N} + o(q^{2N}),
\end{equation}
where
$\zeta_{\Fqt}(s)$ is the Dedekind zeta function of $\Fqt$, or equivalently the zeta function of the algebraic curve
$\PP^1_{\mathbb{F}_q}$.

Datskovsky and Wright's proof uses the adelic Shintani zeta function; in \cite{Z} Zhao obtains another proof
of \eqref{eqn_DW_ff} using algebraic geometry.
A {\itshape trigonal curve} is a 3-fold cover of $\mathbb{P}^1$,
and there is a bijection between isomorphism classes of smooth trigonal curves and cubic rational function fields. Already this 
perspective is enlightening; for example, the Riemann-Hurwitz formula implies that the discriminant
has even
degree, explaining why \eqref{eqn_DW_ff} is an equation for $S(2N)$ and not $S(N)$.
 
The problem, then, is to count smooth trigonal curves. These may be counted by embedding the
curves into {\itshape Hirzebruch surfaces} $F_k$; for each trigonal curve $C$
there is a unique nonnegative integer $k$, called the {\itshape Maroni invariant} of $C$, such that 
$C$ embeds into $F_k = \mathbb{P}_{\mathbb{P}^1}(\calO \oplus \calO(-k)).$

Zhao now counts curves on each surface $F_k$ and sieves for smoothness, using a sieve method related to one developed 
by Poonen \cite{poo}. However, the main technique is different and new. 
The main term in \eqref{eqn_DW_ff} comes from estimating the number of trigonal curves $C$ with
each Maroni invariant $k$, and then summing the results over all integers $k$.
However, the Maroni invariant of a smooth trigonal curve $C$ is at most $N/3$, provided that
the field corresponding to $C$ has absolute discriminant $q^{2N}$. Limiting our sum to $k \leq N/3$ 
naturally introduces a secondary term of order $q^{5N/3}$ into \eqref{eqn_DW_ff}!

Zhao's method, like all other methods described in this paper, comes with error terms; at present
he has obtained an error term of $O(q^{7N/4})$ in \eqref{eqn_DW_ff}. He has some optimism that this error term
can be reduced
below $q^{5N/3}$, thereby proving the secondary term.

\section{Conclusion}
Cubic fields have seen a great deal of recent study recently, and we recommend the papers
of  Cohen-Morra \cite{CM} and Martin-Pollack \cite{MP} (among others) 
for recent studies of cubic fields from other perspectives.
 
Many interesting open questions remain in the subject. For example, what is the ``correct'' error term in the 
Davenport-Heilbronn theorems? Even if we cannot prove it, it would still be interesting
to determine the expected order of magnitude. To get the best error terms it seems that we should count fields
$K$ of degree {\itshape at most} 3, each weighted by $\frac{1}{|\Aut(K)|}$. The data in \cite{R} suggests that the true
error may be smaller than $X^{1/2}$, and a comparison with the divisor problem suggests that the error might be
on the order of $X^{3/8}$. Nevertheless, the comparison with the divisor problem is not exact, and the numerical
data is inconclusive.

Perhaps the most compelling open question concerns quartic
and quintic fields. Asmyptotic formulas were proved by Bhargava \cite{B_quartic, B_quintic};
should these formulas have secondary terms? It is generally believed that they should;
however, so far we lack even a good conjecture. The zeta function approach seems likely to work, 
and Yukie \cite{Y} has made some progress in this direction; the 
geometric approach seems likely to work as well. However, for now the difficulties with either approach
appear rather severe.

Another open question concerns the {\itshape multiplicity} of cubic fields. It is believed that
there should be $\ll n^{\epsilon}$ cubic fields of discriminant $\pm n$, but the best
bound known is $O(n^{1/3 + \epsilon})$, due to Ellenberg and Venkatesh \cite{EV}. Nontrivial
bounds were also obtained by Helfgott and Venkatesh \cite{HV} and Pierce \cite{pierce}, using a
variety of methods. Any of the methods described in the present paper could potentially yield improvements, but none
of them have succeeded to date. (The present author has attempted improvements by means of Shintani
zeta functions, which have led to a number of instructve and interesting failures.)

Finally, this paper begs the question of whether unexpected connections might be found between the four
perspectives presented here, and some promising preliminary results have been obtained in this direction.
We anticipate that this and related questions will be addressed in the near future.

\section{Acknowledgments}
I would like to thank the many people who have shared their insights on these secondary terms, especially
Manjul Bhargava,
Jordan Ellenberg,
Bob Hough,
Arul Shankar,
Kannan Soundararajan,
Takashi Taniguchi,
Melanie Matchett Wood,
and Yongqiang Zhao. 
I would further like to thank Bhargava, Shankar, Taniguchi, and Zhao for comments on this paper in particular.
Most of all I would like to thank Akshay Venkatesh who originally 
brought this topic to my attention.

In addition, I would like to thank the organizers of Integers 2011 for hosting an outstanding
conference.


\begin{thebibliography}{99}

\bibitem{baily} A. M. Baily, 
\emph{On the density of discriminants of quartic fields,} J. Reine Angew. Math. 
\textbf{315} (1980), 190--210. 



\bibitem{BBP} K. Belabas, M. Bhargava, and C. Pomerance,
\emph{Error estimates in the Davenport-Heilbronn theorems}, Duke Math. J. \textbf{153} (2010), 
no. 1, 173--210.


\bibitem{B_quartic} M. Bhargava, 
\emph{The density of discriminants of quartic rings and fields},
Ann. of Math. (2) \textbf{162} (2005), no. 2, 1031--1063. 

\bibitem{B_deg_n} M. Bhargava,
\emph{
Mass formulae for extensions of local fields, and conjectures on the density of number field discriminants},
Int. Math. Res. Not. (2007), no. 17, Art. ID rnm052, 20 pp. 

\bibitem{B_icm} M. Bhargava,
\emph{
Higher composition laws and applications},
Proceedings of the International Congress of Mathematicians, Vol. II, 271Ð294, 
Eur. Math. Soc., Z\"urich, 2006.

\bibitem{B_quintic} M. Bhargava,
\emph{The density of discriminants of quintic
rings and fields}, Ann. Math. (2) \textbf{172} (2010), 
no. 3, 1559-1591.

\bibitem{BST} M. Bhargava, A. Shankar, and J. Tsimerman,
\emph{On the Davenport-Heilbronn theorem and second order terms}, preprint; available at
\url{http://arxiv.org/abs/1005.0672}.

\bibitem{BTT} M. Bhargava, T. Taniguchi, and F. Thorne, 
work in preparation.


\bibitem{CDO_quartic} 
H. Cohen, F. Diaz y Diaz, and M. Olivier, 
\emph{Enumerating quartic dihedral extensions of $\Q$}, 
Compositio Math. \textbf{133} (2002), 65--93. 

\bibitem{CDO} H. Cohen, F. Diaz y Diaz, and M. Olivier,
\emph{Counting discriminants of number fields},
J. Th\'eor. Nombres Bordeaux \textbf{18} (2006), no. 3, 573--593.

\bibitem{CL} H. Cohen and H. Lenstra,
\emph{Heuristics on class groups of number fields},
Number theory, Noordwijkerhout 1983, 33--62, 
Lecture Notes in Math., 1068, Springer, Berlin, 1984.

\bibitem{CM} H. Cohen and A. Morra,
\emph{Counting cubic extensions with given quadratic resolvent}, J. Algebra \textbf{325} (2011), 461--478. 



\bibitem{DW2} B. Datskovsky and D. Wright,
\emph{The adelic zeta function associated to the space of binary cubic forms. II. Local theory},
J. Reine Angew. Math. \textbf{367}  (1986), 27--75.

\bibitem{DW3} B. Datskovsky and D. Wright,
\emph{Density of discriminants of cubic extensions},
J. Reine Angew. Math. \textbf{386}  (1988), 116--138. 


\bibitem{DH} H. Davenport and H. Heilbronn,
\emph{On the density of discriminants of cubic fields. II}, 
Proc. Roy. Soc. London Ser. A \textbf{322}  (1971), no. 1551, 405--420. 

\bibitem{DF} B. N. Delone and D. K. Faddeev,
\emph{The theory of irrationalities of the third degree (in English translation)},
AMS, Providence, 1964.


\bibitem{duke} 
W. Duke,
\emph{Hyperbolic distribution problems and half-integral weight Maass forms},
Invent. Math., \textbf{92} (1):73--90, 1988.

\bibitem{EV_fields} J. Ellenberg and A. Venkatesh,
\emph{The number of extensions of a number field with fixed degree and bounded discriminant},
Ann. of Math. \textbf{163} (2006), no. 2, 723--741.

\bibitem{EV} J. Ellenberg and A. Venkatesh,
\emph{Reflection principles and bounds for class group torsion},
Int. Math. Res. Not. no. 1 (2007), Art. ID rnm002.

\bibitem{GGS} W. T. Gan, B. Gross, and G. Savin,
\emph{Fourier coefficients of modular forms on $G_2$},
Duke Math. J. \textbf{115} (2002), 105--169.

\bibitem{GS} E. Golod and I. Shafarevich,
\emph{On the class field tower} (Russian),
Izv. Akad. Nauk SSSR Ser. Mat. \textbf{28} (1964), 261--272.

\bibitem{HV} H. Helfgott and A. Venkatesh,
\emph{
Integral points on elliptic curves and 3-torsion in class groups}, 
J. Amer. Math. Soc. \textbf{19} (2006), no. 3, 527--550.

\bibitem{H} B. Hough,
\emph{Equidistribution of Heegner points associated to the 3-part of the class group},
preprint; available at
\url{http://arxiv.org/abs/1005.1458}.


\bibitem{L}
F. Lemmermeyer,
\emph{Class field towers},
\url{http://www.rzuser.uni-heidelberg.de/~hb3/publ/pcft.pdf}.

\bibitem{MP}
G. Martin and P. Pollack,
\emph{
The average least character nonresidue and further variations on a theme of Erd\H{o}s},
preprint; available at \url{http://arxiv.org/abs/1112.1175}.

\bibitem{M}
J. Martinet, 
\emph{
Petits discriminants}, 
Ann. Inst. Fourier \textbf{29} (1979), 159--170.



\bibitem{N} J. Neukirch,
\emph{Algebraic number theory},
Springer-Verlag, Berlin, 1999.

 
\bibitem{O} A. Odlyzko,
\emph{
Bounds for discriminants and related estimates for class numbers, regulators and zeros of zeta functions: 
a survey of recent results.}, S\'em. Th\'eor. Nombres Bordeaux (2) \textbf{2} (1990), no. 1, 119--141.

\bibitem{pierce}
L. Pierce,
\emph{The 3-part of class numbers of quadratic fields},
J. London Math. Soc. (2) \textbf{71} (2005), no. 3, 579--598. 

\bibitem{poo}
B. Poonen,
\emph{Bertini theorems over finite fields},
Ann. Math. (2) \textbf{160} (2004), no. 3, 1099Ð1127. 

\bibitem{R} D. Roberts,
\emph{Density of cubic field discriminants},
Math. Comp. \textbf{70}  (2001),  no. 236, 1699--1705.


\bibitem{SK} M. Sato and T. Kimura,
\emph{A classification of irreducible prehomogeneous vector spaces and their relative invariants}.
Nagoya Math. J. \textbf{65} (1977), 1--155. 

\bibitem{SS} M. Sato and T. Shintani,
\emph{On zeta functions associated with prehomogeneous vector spaces},
Ann. of Math. (2) \textbf{100}  (1974), 131--170. 

\bibitem{S} T. Shintani,
\emph{On Dirichlet series whose coefficients are class numbers of integral binary cubic forms},
J. Math. Soc. Japan  \textbf{24} (1972), 132--188. 


\bibitem{sound} K. Soundararajan,
\emph{Divisibility of class numbers of imaginary quadratic fields},
J. London Math. Soc. (2) \textbf{61} (2000), 681--690.


\bibitem{TT} T. Taniguchi and F. Thorne,
\emph{Secondary terms in counting functions for cubic fields}, submitted; available at
\url{http://arxiv.org/abs/1102.2914}.

\bibitem{TT_L}T. Taniguchi and F. Thorne,
\emph{Orbital $L$-functions for the space of binary cubic forms}, submiited; available at
\url{http://arxiv.org/abs/1112.5030}.

\bibitem{tate} J. Tate,
\emph{Fourier analysis in number fields and Hecke's zeta-functions}, thesis (Princeton, 1950); reprinted in
J. W.S. Cassels and A. Fr\"ohlich, eds., 
\emph{Algebraic number theory,}
Academic Press, London, 1986.




\bibitem{wong} 
S. Wong, 
\emph{Automorphic forms on $\GL(2)$ and the rank of class groups,}
J. Reine Angew. 
Math. \textbf{515} (1999), 125--153. 

\bibitem{DW1} D. Wright,
\emph{The adelic zeta function associated to the space of binary cubic forms. I. Global theory},
Math. Ann.  \textbf{270}  (1985),  no. 4, 503--534.

\bibitem{W_iwasawa} D. Wright,
\emph{Twists of the Iwasawa-Tate zeta function}, Math. Z. \textbf{200} (1989), 209--231.

\bibitem{WY} D. Wright and A. Yukie,
\emph{Prehomogeneous vector spaces and field extensions},
Invent. Math. \textbf{110} (1992), no. 2, 283--314.

\bibitem{Y} A. Yukie,
\emph{Shintani zeta functions},
London Mathematical Society Lecture Note Series, 183, Cambridge University Press, Cambridge, 1993.

\bibitem{Z} Y. Zhao,
doctoral thesis, University of Wisconsin, in preparation.

\end{thebibliography}
\end{document}